\documentclass[a4paper,11pt]{article}

\usepackage{amsmath,amsfonts,amssymb,amsthm}
\usepackage[latin1]{inputenc}

\usepackage{enumerate}
\usepackage{latexsym,amssymb,amsthm,amsxtra}
\usepackage{amsmath,mathrsfs}
\usepackage{tikz}
\usepackage{pgf}
\usepackage{dcolumn}
\usepackage{stmaryrd}
\usepackage{pstricks}
\usepackage{pst-node}
\usepackage{multido}

\usepackage[textsize=scriptsize]{todonotes}

\def\ep{\hfill $\Box$}
\def\ie{i.e.\,}
\def\bp{\noindent{\bf Proof.}\ }
\def\epsilon{\varepsilon}

\usepackage{anysize}
\marginsize{1.75cm}{2.25cm}{2.5cm}{2.5cm} 
\newtheorem{theorem}{Theorem}[section]

\newtheorem{prop}[theorem]{Proposition}
\newtheorem{coro}[theorem]{Corollary}
\newtheorem{remark}[theorem]{Remark}

\newtheorem{assumption}[theorem]{Assumption}
\newtheorem{example}[theorem]{Example}

\newcommand{\be}{ \begin{equation}}
\newcommand{\ee}{\end{equation}}

\def\E{{\mathbb E}}

\def\P{{\mathbb P}}

\def\G{{\mathcal{G}}}

\def\({{\Bigl(}}
\def\){{\Bigr)}}

\def\square{\ifmmode\sqr\else{$\sqr$}\fi}
\def\sqr{\vcenter{
         \hrule height.1mm
         \hbox{\vrule width.1mm height2.2mm\kern2.18mm\vrule width.1mm}
         \hrule height.1mm}}                  

\title{Scaling limits for exploration algorithms}

\author{P. Bermolen, M. Jonckheere, J. Sanders}

\begin{document}

\maketitle

\abstract{We consider an exploration algorithm where at each step, a random number of items become active while related items get explored. Given an initial number of items $N$ growing to infinity and building on a strong homogeneity assumption, we study using scaling limits of Markovian processes statistical properties of the proportion of active nodes in time. This is a companion paper that rigorously establishes the claims and heuristics presented in \cite{sanders_sub-poissonian_2015}.}

\section{Introduction}

Assume there exists a binary relation between items $V =\{ 1, \ldots, N \}$, to which we associate a graph where nodes are items such that two items are neighbors if they are related. Let $A_t$ be the set of active items at time $t$ (step $t$) and $B_t$ the set of explored items.
We assume that initially, $A_0=B_0=\{\emptyset\}$. Then, we consider the following exploration process:
(i) select $I^t \subset V \setminus {A_t \cup B_t}$ and determine its neighbors in the set of nonexplored items,
$\mathcal{N}^{I_t} \subset V \setminus {A_t \cup B_t}$, and (ii) actualize $A_t$ and $B_t$ by setting
\begin{gather*}
A_{t+1}= A_t\cup \{I_t\}, \\
B_{t+1}= B_t\cup \mathcal{N}^{I_t}.
\end{gather*}

These exploration algorithms can be used for instance to approximate the evolution of parking processes \cite{bermolen}
and as we shall show, classes of random sequential adsorption processes. In case $I^t$ is a single item at each step, the algorithm 
discovers in a greedy manner independent sets of the relation graph. For instance, this type of exploration algorithm is used for defining subsets of communicating nodes in communication networks with interferences. Otherwise, one can think of $I^t$ as a subgraph, and $B_t$ the set of neighbors of this subgraph. Again, it might be linked to communications procedures where parts of a network might get priority to transmit.

In such problems, the relation graph might be the outcome of spatial effects (nodes interacting though a geometry which can be itself random) and purely random relations between nodes. This is definitely the case of wireless networks which have radio conditions defining the level of admissible interference between two competing nodes. A hardcore interference graph model might hence define an edge between two nodes if their radio conditions would impede a synchronous communication. 

Blockade effects in complex systems of interacting particles can also be described using relation graphs. Of particular interest to us are specially prepared gases that consist of ultracold atoms that can reach a ``Rydberg state'' from a ``ground state''. The essential feature of these particles is that each atom that is in its Rydberg state, prevents neighboring atoms from reaching their Rydberg state. This is similar in spirit to the interference constraints in wireless networks \cite{sanders_wireless_2014}, and the essential features of the blockade effect can therefore again be described using interference graphs. This realization also allows us to study statistical properties of the proportion of atoms ultimately in the Rydberg state through the exploration process defined above \cite{sanders_sub-poissonian_2015}.

  Note that in general, the dimension needed to represent the exploration algorithm as a Markov process is $n$, the size of the graph which impedes simple computations.
In this article, we suppose a strong homogeneity assumption on the relation between items
so that the dynamics can be made Markov in dimension $1$ which is a very crude simplification for most problems but which has the great advantage to lead to simple tractable processes.
We here use classical tools (fluid limits and diffusion approximations) to derive computable characterizations of the performance of such algorithms under such an assumption. 
These results allow us to prove functional laws of large numbers and central limit theorems for the proportion of active items.

This article is structured as follows. In Section \ref{sec:rsa}, we study the case of random sequential adsorption (RSA) under an homogeneity assumption. We first review usual functional law of large numbers and central limit theorems for scaled Markov processes aiming at explicit error bounds. We then use these results to study hitting times that capture statistical properties of the proportion of active nodes.

\section{Random adsorption under an homogeneous relation}\label{sec:rsa}


Assume from here on that precisely one item is selected in each step. Assume also the following homogeneity assumption on the graph $\G$ induced by the binary relation between items.
\begin{assumption}
If $(\G_1, \G_2)$ is a partition of $\G$, $\forall$ $i \in \G_2$,
the distribution of the number of neighbors of $i$ in $\G_2$ depends only
on $|\G_1|$ and $|\G_2|$.
\label{A1}
\end{assumption}

\begin{remark}
Although Assumption (\ref{A1}) is not valid in many practical cases such as random geometric graphs and random graphs with generic degree distribution, Assumption (\ref{A1}) is crucial to get a one-dimensional analysis. It can however be considered a reasonable approximation for many systems, and it is satisfied for instance by Erd\"os--R\'enyi's random graph. See \cite{bermolen} for a study of scaling limits in infinite dimension which applies to a much larger class of problems.
\end{remark}

Let $Z_n$ be the number of explored items at step $n$, \ie $Z_n= \sharp \{A_n\cup B_n\}$. Then,
\begin{gather*}
Z_n = Z_{n-1}+1+\xi_n \quad \text{and} \quad Z_0=0,
\end{gather*}
where $\xi_n$ is the distribution of the number of neighbors that an item has at step $n$ in the remaining non-explored portion of the graph.

Under Assumption (\ref{A1}) the distribution of $\xi_n$ depends only on $Z_{n-1}$, which with a slight abuse of notation we denote by $\xi_{Z_{n-1}}$. This also implies that $Z=\{Z_n\}_{n\in\mathbb{N}}$ is a discrete Markov process taking values in $\{0,\dots, N\}$, and that $Z$ is an increasing process with in $N$ an absorbing state. The transition probabilities are given by
\begin{equation*}
p_{xy}(n)=P(Z_n=y|Z_{n-1}=x)=P(\xi_x=y-x-1)\quad \text{with} \quad y > x.
\end{equation*}
If we now denote by $(p_N(\cdot,x))$ the distribution of the number of neighbors in $\G_2$ of a given vertex $i\in \G_2$ given that $(\G_1,\G_2)$ is a partition of $\G$ with $|\G_1|=x$, the transition probabilities can be written as
\begin{gather*}
p_{x, x+k+1}=P(\xi_x=k)=p_N(k,x) \quad \text{with}\quad  k\geq 0.
\end{gather*}
The transition probabilities in case of the Erd\"os--R\'enyi random graph are given by the Binomial distribution, \ie $p_N(k,x) = \binom{N-x-1}{k} p^k (1-p)^{N-x-k-1}$.

\subsection{Preliminaries}
\subsubsection{ Functional law of large numbers}

Given a partition $(\G_1,\G_2)$ of $\G$ such that $|\G_1|=x$, we consider the mean $\gamma_N(x)$ and variance $\psi_N(x)$ of the number of neighbors in $\G_2$ of a given vertex $i\in\G_2$
\begin{gather*}
\gamma_N(x)=\sum_{k=0}^{N-1} k p_N(k,x),\\
\psi_N(x)=\sum_{k=0}^{N-1} (k-\gamma_N(x))^2 p_N(k,x),
\end{gather*}
and define $\bar{\gamma}_N= \sup_{x} \gamma_N(x)$, $\bar{\psi}_N = \sup_{x} \psi_N(x).$

We now consider the scaled process on time and space.
Define the scaled process, viewed as a piece-wise constant trajectory process in continuous time, i.e.,
for all $t>0$
\begin{gather*}
Z^{N}_t
= \frac{ Z_{[tN]} }{N}
\end{gather*}
Here, $[x]$ is the integer part of $x$, and we suppose that $Z_0=0$.
We will now derive a law of large numbers for $Z_t^N$ using classical tools \cite{Darling2008}. The proof of convergence relies on classical techniques, which we leverage to obtain error bounds along the way.

\begin{prop}
\label{prop:LLN_comparing_Znt_and_Zt}
\label{prop:fluid_limit_disc}
If there exists a ($C_L$)-Lipschitz function $\gamma$ on $\mathbb R^+$ such that
\begin{gather}
\sup_{x} \left|\gamma_N(x)-\gamma\(\frac{x}{N}\) \right| \le \delta_N, 
\label{eq:def_gamma}
\end{gather}
then for $p > 1$,
\begin{equation}
\lVert \sup_{s \in [0,T]}|Z^{N}_s- z(s)| \ \rVert_p \le\exp\{ C_L T \} \big(  \delta_N T +   \kappa_p \lVert M_T^N \rVert_p    \big),
\label{eqn:Bound_on_Lp_norm_of_supremum_distance_ZsN_and_zs}
\end{equation}
where $\kappa_p = p / (p-1)$. In particular for $p = 2$,
\begin{gather*}
\lVert \sup_{t \in [0,T]}|Z^{N}(t)- z(t)| \, \rVert_2 \leq \omega_N,
\end{gather*}
where $\omega_N =(\delta_N T+2\sqrt{2 C_N T})\exp\{C_L T\}$, $C_N = \bar{\psi}_N / N$, and $z(t)$ is the solution of the deterministic differential equation
\begin{gather}
\dot{z}= 1+\gamma(z).
\label{eq:ODE.disc}
\end{gather}
\end{prop}

\bp
Using the martingale decomposition of the Markov process $Z$, we have
\begin{gather} \label{eq:martingale_decomp_Z}
Z_l= Z_0+ \sum_{i=0}^l \bigl( 1 + \gamma_N(Z_i) \bigr) + M_l,
\end{gather}
where $M_l$ is a local martingale which is actually a global martingale since the state space is finite. 

Scaling \eqref{eq:martingale_decomp_Z}, and viewing its trajectory as piece-wise constant, it follows that
\begin{align}
Z^N_t
= \frac{Z_{[tN]}}{N}
&=  \nonumber Z^N_0+ \frac{1}{N} \sum\limits_{i=0}^{[tN]} \bigl( 1+\gamma_N(Z_i) \bigr) + \frac{M_{[tN]}}{N} \\
&= \nonumber  Z^N_0+ \frac{1}{N} \int_{0}^{[tN]} \bigl( 1+\gamma_N(Z_s) \bigr) ds + \frac{M_{[tN]}}{N}\\
&= Z^N_0 + \int_{0}^{[t]} \bigl( 1 + \gamma_N(Z_{uN}) \bigr) du + M^N_t,
\label{eq:martingale_decomp_Zn}
\end{align}
where the latter equality follows from a change of variables, and an introduction of notation for the scaled martingale,
$
M^N_t
= M_{[tN]} / N
$.


Using the integral version of \eqref{eq:ODE.disc}, the triangle inequality, and Lipschitz continuity of $\gamma$, we find that
\begin{align}
&
\sup_{s\in[0,t]}| Z^N_s - z(s) |
\leq
\sup_{s \in [0,t]} \Bigl( | Z^N_0 - z(0) | + \int_0^{[t]} \bigl| \gamma_N(Z_{uN}) - \gamma(z(u)) \bigr| du + | M_t^N | \Bigr)
\nonumber \\ &
\leq |Z_0^N-z(0)| + C_L \int_0^{[t]} \sup_{u\in[0,s]}| Z^N_u - z(u)| ds + \delta_N t + \sup_{s\in[0,t]}|M_s^N| 
\nonumber \\ &
\leq |Z_0^N-z(0)| + C_L \int_0^{t} \sup_{u\in[0,s]}| Z^N_u - z(u)| ds + \delta_N t + \sup_{s\in[0,t]}|M_s^N|.
\label{eq:epsilon_n}
\end{align}
Define $\epsilon_N(T)= \underset{s \in [0,T]}{\sup}|Z^{N}_s- z(s)|$, so that from \eqref{eq:epsilon_n} it follows that
\begin{equation*}
\epsilon_N(T) \leq |Z_0^N-z(0)| + \delta_N T + \sup_{s\in[0,t]}|M_s^N| + C_L \int_0^T \epsilon_N(s) ds.
\end{equation*}
Recall that $Z_0^N = z(0) = 0$, and because $\delta_N T + \sup_{s\in[0,T]}|M_s^N|$ is nondecreasing in $T$, it follows from Gr\"{o}nwall's lemma that
\begin{gather*}
\epsilon_N(T)
\leq \bigl( \delta_N T + \sup_{s\in[0,T]}|M_s^N| \bigr) \exp\{ C_L T\}.
\end{gather*}
Using Minkowsky's inequality for $p \in [1,\infty)$, strict monotonicity of $\exp\{ C_L T \}$ and $\delta_N T$, and the triangle inequality, we find
$$ \lVert \epsilon_N(T) \rVert_p \le\exp\{ C_L T \} \bigl( \delta_N T + \lVert \sup_{s \in [0,T]} |M_s^N| \rVert_p \bigr).$$
Finally, using Doob's martingale inequality for $p > 1$, we obtain
\begin{equation}
\lVert \epsilon_N(T) \rVert_p \le\exp\{ C_L T \} \big( \delta_N T + \kappa_p \lVert M_T^N \rVert_p \big).
\end{equation}

In $L^2$, this inequality can be further simplified by computing the increasing process associated to the martingale, \ie for $l \geq 0$,
$$
\E[ (M_l)^2 ] 
= \E[ \langle M_l\rangle ]
= \E\Bigl[ \sum_{i=0}^{l} {\rm Var}[ \gamma_N(Z_i) ] \Bigr].
$$
where, 
\begin{gather*}
{\rm Var}[ \gamma_N(x) ] 
= \sum_{k=0}^{N-x-1} (k+1)^2 p_{x, x+k+1} - \Bigl( \sum_{k=0}^{N-x-1} (k+1) p_{x, x+k+1} \Bigr)^2 
= \psi_N(x).
\end{gather*}
Then,
\begin{equation}
\lVert M_t^N \rVert_2^{2} 
= \E[(M_{t}^N)^2] 
= \frac{\E[M^2_{ [ tN ] }]}{N^2}=\frac{1}{N^2}\sum_{i=0}^{[tN]} \psi_N(Z_i)\leq C_N t.
\label{eqn:Bound_on_L2_norm_of_MtN}
\end{equation}
This completes the proof.
\ep

\begin{coro}\label{coro.fluid.disc}
If the distribution of the number of neighbors is such that $\delta_N \to 0$ as $N \to \infty$ and $\bar{\psi}_N = o(N)$, the scaled process $Z^N_t$ converges to $z_t$ in $L^1$ uniformly on compact time intervals.
\end{coro}

\begin{coro}\label{coro.fluid.disc.Nrandom}
If the number of initial items $N$ is itself random and independent of the trajectory of $Z$, meaning that $Z$ can be constructed (i) as a functional of $N$ and (ii) of other random variables that are independent of $N$, then
\begin{gather*}
\E \bigl[ \sup_{t \in [0,T]}|Z^{N}(t)- z(t)| \bigr] 
\leq \E(\delta_N T) +2 \E \sqrt{{ C_N T}})\exp\{ C_L T\}.
\end{gather*}
\end{coro}

\begin{example}[Sparse Erd\"{o}s--R\'enyi Graph with a Poissonian number of vertices]\label{ex:erdos_renyi_disc}
Suppose that given $N$, the graph $\G=\G(N, c/N)$ is a sparse Erd\"{o}s--R\'enyi graph, \ie
$p_N(\cdot,x)$ is the probability mass function of the binomial distribution ${\rm }Bin(N-x-1, c/N)$ with $c > 0$.
Additionally, suppose that $N-1$ is Poisson distributed with parameter $h$. The mean and variance of $p_N(.,x)$ are then given by
\begin{equation}
\gamma_N(x)= (N-x-1) \frac{c}{N}, \quad
\psi_N(x)= (N-x-1) \frac{c}{N} \Bigl( 1- \frac{c}{N} \Bigr). \nonumber
\end{equation}

Define $\gamma(x)=c(1-x)$. Condition \eqref{eq:def_gamma} is satisfied and as Lipschitz constant $C_L=c$ suffices. Moreover, $\delta_N= {c / N}$ and $\bar{\psi}_N \le c$. In this case the deterministic differential equation in \eqref{eq:ODE.disc} reads
\begin{gather*}
\dot{z} = 1+ c(1-z)=(1+c) - cz,
\end{gather*}
which can be explicitly solved, resulting in
\begin{gather*}
z(t)=\rho+(z_0 -\rho)e^{- t}\quad \mbox{with} \quad \rho=\frac{1+c}{c}.
\end{gather*}
Since $z(0)=0$, $z(t)=\rho(1 -e^{- t})$. Observe that $\underset{t\to\infty}{\lim} z(t)=\rho>1$.

Furthermore, from Corollary \ref{coro.fluid.disc.Nrandom} we obtain using Cauchy-Schwarz's inequality that there exists a constant $C_1$ such that
\begin{align*}
\E \Bigl[ \sup_{t \in [0,T]}|Z^{N}(t)- z(t)| \Bigr]
& 
\le  c \E(1/N)  +   \E(1/\sqrt{N})   2 \sqrt{ c T}\exp\{ c T\}
\\ &
\le  C_1 \E(1/N)^{1/2}  
= C_1 \Big(  \frac{1-\exp(-h)}{h} \Big)^{1/2}.
\end{align*}

\end{example}

%

\subsubsection{Diffusion approximations with errors bounds}

We now proceed and derive a functional central limit theorem. The convergence proof again uses classical techniques, which we use to determine error bounds. To that end, we apply results of \cite{kurtz78} which are based on results by Koml\'{o}s-Major-Tusn\'{a}dy. These results allow one to construct a Brownian motion and either a Poisson process or random walk on the same probability space. Since we are concerned with discrete time, we need to consider the random walk case, see also \cite{berkes2014}. In order to obtain explicit error bounds, we suppose stronger assumptions on the transitions probabilities than would be needed when only proving convergence.

\begin{prop}
\label{prop:Diffusion_approximation_with_error_bound}
If there exists a function $p$ on $(\mathbb N,\mathbb R_+)$, and a sequence $(\epsilon_k)_{k=0,1,\ldots}$ such that
\begin{gather*}
|p_N(k,[Nx])- p(k,x)| \le \frac{\epsilon_k}{N}, \\
|p(x,x+k)- p(y,y+k)| \le M \epsilon_k |x-y|, \\
\sum_{k} k^2 |p(x,x+k)^{1/2}- p(y,y+k)^{1/2}|^2 \le M  |x-y|^2, \\
\sum_{k} k \epsilon_k^{1/2} < \infty,
\end{gather*}
and if $\gamma$ is twice differentiable with bounded first and second derivatives,
then the process 
\begin{equation}
W^N_t=\sqrt{N}\(Z^N_t-z(t)\)
\end{equation}
converges in distribution towards $W_t$, the unique solution
of the stochastic differential equation
\begin{gather}
dW(t)= \gamma'(z(t))W(t)dt+  \sqrt{\beta'(t)} dB_1(t).
\label{eq:STDE}
\end{gather}
Here, $B_1(t)$ denotes a standard Brownian motion, $\beta(t)=  \int_{0}^t  \psi(z(s)ds$, and $z(t)$ is the solution of \eqref{eq:ODE.disc}. 
Furthermore, 
\begin{equation}
\E\(\sup_{t \le T} |W^N_t-W_t|\) \le C \frac{\log(N)}{\sqrt{N}}.
\end{equation}

\end{prop}

\bp
We adapt the results of Kurtz which were derived for continuous time Markov jump processes.
For doing so, we can
replace the Poisson processes  involved in the construction of the jump processes by some random walks
that can be used to construct discrete time Markov chains.
We can then use exactly the same steps as in \cite{kurtz78},
by first comparing the original process $Z^N$ with a diffusion of the form
$$\tilde Z^N_t = \frac{1}{N} \sum_{l \le N} l B_l (N \sum_0^t p_N(l,\tilde Z^N_s)ds),$$
that is a sum of a finite number of scaled independent Brownian motions $B_l$.

Rewriting the inequalities in \cite[(3.6)]{kurtz78}, and using a random walk version of the approximation lemma of  Koml\'os-Major-Tusn\'ady \cite{berkes2014}, we obtain
\begin{equation}
\E \(\sup_{t \le T} |\tilde Z^N_t - Z^N_t| \)\le C_2 \frac{\log(N)}{N}.
\end{equation}
This leads using the results of \cite[Section~3]{kurtz78} to
\begin{equation}
\E \(\sup_{t \le T} | W^N_t- W_t | \)\le C_3 \frac{\log(N)}{\sqrt{N}},
\end{equation}
which concludes the proof.
\ep

\begin{example}[ER case -- Continued]
The relation between the binomial coefficients and the Poisson distribution are well studied.
Defining
$$p_N(k,[xN])= \binom{N-[Nx]-1}{k} \Bigl( \frac{c}{N} \Bigr)^k \Bigl( 1-\frac {c}{N} \Bigr)^{N-[Nx]-1},$$
and using (for instance) the Stein-Chen method \cite{teerapabolarn}, we have that
$$ |p_N(k,[xN])-p(k,x)| \le \frac{c}{N} p(k,x) ,$$
which shows that the assumptions of Proposition~\ref{prop:Diffusion_approximation_with_error_bound} are satisfied.
\end{example}

\subsection{LLN and CLT for the Hitting time}

If $N<\infty$, the exploration algorithm finishes at
\begin{equation}
T_N^* = \inf\{ \tau \in \mathbb{N}_+ | Z_\tau = N \} \leq N < \infty.
\end{equation}
This time $T_N^*$ is a hitting time for the Markov process. Since the algorithm adds precisely one node at each step, we have that the final number of active items is exactly $T^*_N$, \ie $A_{T_N^*} = T_N^*$. Because we wish to determine the statistical properties of $A_{T_N^*}$, we will seek not only a first-order approximation for $T^*_N$, but also prove a central limit theorem result as the initial number of items $N$ goes to infinity. From here on onward, we denote by $T^*$ the solution to $z(T^*)=1$.

\begin{prop}
\label{prop:Bound_on_L2_norm_of_TNstar_minus_Tstar}
For all $\delta > 0$, there exists a constant $C_\delta$ depending only on $\delta$, so that
\begin{equation}
\P\Bigl[ \Bigl| \frac{T^*_N}{N} - T^* \Bigr| \ge \delta \Bigr]
\le C_\delta \omega_N.
\end{equation}
Furthermore, if there exist constants $\varepsilon, c_1 > 0$ so that $\gamma(z(s)) \le  1- \epsilon$ for all $s \ge c_1$,
then there exists a constant $C$ such that
\begin{gather*}
\Big\lVert \frac{T^*_N}{N} - T^* \Big\rVert_2
\leq C \omega_N. 
\end{gather*}

\end{prop}

\bp
Remark that if $| z(s) - Z_{s}^N | \le \delta$, then there exists a finite constant $A>0$ such that
\begin{equation}
\Bigl|\frac{T^*_N}{N} - T^* \Bigr| 
\le \sup_{ t \in [z^{-1}(1) - \delta, z^{-1}(1) + \delta] } z^{-1}(t) 
= \sup_{t \in [T^*-\delta, T^*+\delta]} z^{-1}(t) 
\le A \delta.
\end{equation}
Hence the first claim follows directly from the observation that the event
\begin{equation}
\Big\{ \Bigl| \frac{T^*_N}{N} - T^* \Bigr| \ge A \delta \Big\} \subset \big\{ | z(s) - Z_{s}^N | \ge  \delta \big\}.
\end{equation}

Now since $Z_0=0$, and using that $Z_{T_N^*} / N = z(T^*) = 1$ together with \eqref{eqn:Bound_on_Lp_norm_of_supremum_distance_ZsN_and_zs} and \eqref{eq:martingale_decomp_Z}, we find
%
\begin{align}
\label{eq:Tmart}
&
\frac{T^*_N}{N}- T^*
= \int_0^{T^*} \gamma(z(s)) ds - \int_0^{\frac{T^*_N}{N}} \gamma_N(Z_{sN})ds 
- \frac{M_{T_N^*}}{N}
\\ &
= \int_0^{ \frac{T_N^*}{N} \wedge T^* } ( \gamma(z(s)) - \gamma_N(Z_{sN}) ) ds + \int_{ \frac{T_N^*}{N} \wedge T^* }^{  T^* } \gamma(z(s)) ds - \int_{ \frac{T_N^*}{N} \wedge T^* }^{  \frac{T_N^*}{N} } \gamma_N(Z_{sN}) ds 
- \frac{M_{T_N^*}}{N}. 
\nonumber
\end{align}
Then, using the triangle inequality,
\begin{align*}
\Bigl| \frac{T^*_N}{N} - T^* \Bigr|
\leq & \int_0^{ \frac{T_N^*}{N} \wedge T^*} |\gamma(z(s)) - \gamma_N(Z_{sN})| ds + \int_{ \frac{T_N^*}{N} \wedge T^*}^{ T^* } | \gamma(z(s))| ds 
\nonumber \\ &
+ \int_{ \frac{T_N^*}{N} \wedge T^* }^{ \frac{T_N^*}{N} } | \gamma_N(Z_{sN}) | ds + |M_{T_N^*}^N|.
\end{align*}
Approximating $\gamma_N$ by $\gamma$ and using the Lipschitz continuity of $\gamma$  and $\max(T_N^* / N, T^*)  \leq 1$,
\begin{align*}
\Bigl| \frac{T^*_N}{N} - T^* \Bigr| &
\leq 2 C_L \sup_{s \le1 } | z(s) - Z_{s}^N | +  \delta_N +  \int_{ \frac{T_N^*}{N} \wedge T^*}^{T^*} | \gamma(z(s)) | ds + \int_{ \frac{T_N^*}{N} \wedge T^* }^{ \frac{T_N^*}{N} } | \gamma(z(s)) | ds   + |M_{T_N^*}^N|.
\end{align*}
Splitting cases, using that $\gamma(z(s)) \leq 1 - \varepsilon$ for $0 < c_1 \leq s $,  the bound $\lVert T^*_N / N - T^* \rVert_1 \leq 1$ and $T_N^* / N  \leq 1$,
\begin{align*}
\Bigl| \frac{T^*_N}{N} - T^* \Bigr|&
\leq  \delta_N + 2 C_L \sup_{ s \le  1 } | z(s) - Z_{s}^N | + (1-\epsilon) \Bigl|\frac{T^*_N}{N} - T^* \Bigr| + C_2  1_{ \frac{T^*_N}{N}  \le c_1 } + |M_{T_N^*}^N|.
]\end{align*}

Now note that similarly to our previous argumentation before, there exists a $\delta$ such that
\begin{equation}\label{eq:pr1}
\Bigl\{ \frac{T^*_N}{N} \le c_1  < T^* \Bigr\}
\subset \{ | z(s) - Z_{s}^N | > \delta \}.
\end{equation} 
Hence, using \eqref{eq:pr1} together with Markov's inequality and the Minkowski inequality, there exists a constant $C_3$ such that
\begin{gather*}
\epsilon \Bigl\lVert \frac{T^*_N}{N} - T^* \Big\rVert_2 \leq  \delta_N + C_3 \lVert \sup_{s \le 1 } | z(s) - Z_{s}^N | \rVert_2 + \lVert M_{T_N^*}^N \rVert_2.
\end{gather*} 
%
Using Proposition~\ref{prop:LLN_comparing_Znt_and_Zt}, \eqref{eqn:Bound_on_L2_norm_of_MtN}, and the fact that $T^*_N/N \leq 1$, we obtain that
\begin{align*} 
\Big\lVert \frac{T^*_N}{N} - T^* \Big\rVert_2 
&
\leq   \delta_N + C_4 || \underset{s\in [0, 1 ]}{\sup} | z(s) - Z_{s}^N| \, ||_2 + ||M_{T_N^*}^N ||_2, 
\nonumber \\ &
\le  \delta_N + C_4 \omega_N + \bar \psi_N/N E(T_N^*/N) \leq C_5 \omega_N.
\end{align*}
This concludes the proof.
\ep

\begin{coro}
Suppose $\gamma(1) \neq 1$.
The random variable $\sqrt{N} ( T_N^* / N - T^* )$ converges in $L^1$ to $W_{T^*}$, and
\begin{gather*}
\E \Bigl[ \Bigl| \sqrt{N} \Bigl( \frac{T_N^*}{N} - T^* \Bigr) (1-\gamma(1))+ W_{T^*} \Bigr| \Bigr] \le \sqrt{N} \omega_N^2.
\end{gather*}
Here, $W_{T^*}$ is a centered Gaussian random variable with variance
\begin{equation}
\sigma^2{ = \frac{m_{T^*}}{ 1 - \gamma(1) }, }
\end{equation}
where $m_t = \E[ W_t^2 ]$ solves the differential system
\begin{gather}
\dot{m_t}=-2\dot{\gamma}(z_t) m_t+ \dot{\beta}(t),
\quad m_0 = 0.
\label{eq:ode_var_er_disc}
\end{gather}
\end{coro}

\bp
First, expand
\begin{align}
&
\Bigl| \sqrt{N} \Bigl( \frac{T_N^*}{N}-T^* \Bigr)( 1-\gamma(1) ) + W_{T^*} \Bigr|
= \Bigl| \sqrt{N} \Bigl( \frac{T_N^*}{N}-T^* \Bigr)( 1-\gamma(1) ) + W_{T^*}^N + W_{T^*} - W_{T^*}^N \Bigr|
\nonumber \\ &
= \Bigl| \sqrt{N} \Bigl( \frac{T_N^*}{N}-T^* \Bigr)( 1-\gamma(1) ) + \sqrt{N} ( Z_{T^*}^N - z(T^*) ) + W_{T^*} - W_{T^*}^N \Bigr|
\nonumber \\ &
= \Bigl| \sqrt{N} \Bigl( \frac{T_N^*}{N}-T^* \Bigr)( 1-\gamma(1) ) + \sqrt{N} \Bigl( \int_0^{T^*} ds - z(T^*) + \int_0^{T^*} \gamma_N(Z_{sN}) ds + M_{T^*}^N \Bigr) + W_{T^*} - W_{T^*}^N \Bigr|
\nonumber \\ &
= \Bigl| \sqrt{N} \Bigl( \frac{T_N^*}{N}-T^* \Bigr)( 1-\gamma(1) ) + \sqrt{N} \Bigl( - \int_0^{T^*} \gamma(z(s)) ds + \int_0^{T^*} \gamma_N(Z_{sN}) ds + M_{T^*}^N \Bigr) + W_{T^*} - W_{T^*}^N \Bigr|
\nonumber
\end{align}
Then, use \eqref{eq:Tmart} to simplify (with the notation that $\int_a^b = - \int_b^a$ when $a > b$)
\begin{align}
\ldots 
= & \Bigl| \sqrt{N} \Bigl( - \int_0^{ \frac{T^*_N}{N} } \gamma_N(Z_{sN}) ds - M_{T^*_N}^N \Bigr) - \gamma(1) \sqrt{N} \Bigl( \frac{T_N^*}{N}-T^* \Bigr) 
\nonumber \\ &
+ \sqrt{N} \int_0^{T^*} \gamma_N(Z_{sN}) ds + \sqrt{N} M_{T^*}^N + W_{T^*} - W_{T^*}^N \Bigr| 
\nonumber \\
= & \Bigl| \sqrt{N} \Bigl( \int_{ \frac{T^*_N}{N} }^{ T^* } \gamma_N(Z_{sN}) ds - M_{T^*_N}^N \Bigr) - \gamma(1) \sqrt{N} \Bigl( \frac{T_N^*}{N}-T^* \Bigr) + \sqrt{N} M_{T^*}^N + W_{T^*} - W_{T^*}^N \Bigr| 
\nonumber
\end{align}
Finally, add and substract $\int_{T_N^* / N}^{T^*} \gamma(z(s)) ds$, and use the triangle inequality to arrive at
\begin{align}
&
\left| \sqrt{N}\left(\frac{T_N^*}{N}-T^*\right)(1-\gamma(1))+ W_{T^*}\right| 
\leq
| W_{T^*} - W^N_{T*} | 
+ \sqrt{N} \Bigl|{\int_{T_N^*/N}^{T^*}} ( \gamma_N(Z_{sN}) - \gamma(z(s)) ) ds \Bigr|
\nonumber \\ & 
+ \sqrt{N} \Bigl| ]{\int_{T_N^*/N}^{T^*}} \gamma(z(s)) ds  - \gamma(1) \left(\frac{T_N^*}{N}-T^*\right) \Bigr|
+ {\sqrt{N} | M^N_{T^*} - M^N_{T^*_N} |.}
\end{align}

We now bound each of the terms on the right one by one. The first term can be bounded using Proposition~\ref{prop:Diffusion_approximation_with_error_bound},
\begin{equation}
\E[ W^N_{T*}- W_{T^*} ] 
\le C_0 \frac{\log(N)}{\sqrt{N}}.
\end{equation}

The second term can be bounded using (in sequence) the triangle inequality, \eqref{eq:def_gamma}, Lipschitz continuity of $\gamma$, and extending the integration range, \ie
\begin{align}
&
\E \Bigl[ \Bigl| \sqrt{N} \int_{T^*}^{T_N^*/N} \gamma_N(Z_{sN}) - \gamma(z(s)) ds \Bigr| \Bigr] 
\nonumber \\ &
 \le \sqrt{N} \E \Bigl[ \mathrm{sign}\Bigl( \frac{T_N^*}{N} - T^* \Bigr) \int_{T^*}^{T_N^*/N} | \gamma_N(Z_{sN}) - \gamma(z(s)) | ds \Bigr] 
\nonumber \\ &
 \le \sqrt{N} \E \Bigl[ \mathrm{sign}\Bigl( \frac{T_N^*}{N} - T^* \Bigr) \int_{T^*}^{T_N^*/N} ( \delta_N + | \gamma( Z_s^N ) - \gamma(z(s)) | ) ds \Bigr] 
\nonumber \\ &
 \leq \sqrt{N} \E \Bigl[ \bigl( \delta_N + C_L \sup_{s \leq 1 }| Z_s^N - z(s) | \bigr) \Bigl| \frac{T_N^*}{N} - T^* \Bigr| \Bigr] 
\nonumber \\ &
 = \sqrt{N} \Bigl( \delta_N \Bigl\lVert \frac{T_N^*}{N} - T^* \Bigr\rVert_1 + C_L \Bigl\lVert \sup_{s \leq 1 }| Z_s^N - z(s) | \Bigl( \frac{T_N^*}{N} - T^* \Bigr) \Bigr\rVert_1.
\end{align}
Then using H\"{o}lder's inequality, and the bound $\lVert T^*_N / N - T^* \rVert_1 \leq 1$, we find that
\begin{align}
\ldots & \le \sqrt{N} \Bigl( \delta_N  \Bigl\lVert \frac{T_N^*}{N} - T^* \Bigr\rVert_1  + C_2  \Bigl\lVert \sup_{s \le 1 } | z(s) - Z_{s}^N | \Bigr\rVert_2 \Bigl\lVert \frac{T^*_N}{N} - T^* \Bigr\rVert_2 \Bigr)
\nonumber \\ &
\le \sqrt{N} \Bigl( \delta_N + C_2  \Bigl\lVert \sup_{s \le 1 } | z(s) - Z_{s}^N | \Bigr\rVert_2 \Bigl\lVert \frac{T^*_N}{N} - T^* \Bigr\rVert_2 \Bigr).
\end{align}
%
Hence
\begin{align}
&
\E \Bigl[ \Bigl| \sqrt{N}\int_{T^*}^{T_N^*/N} \gamma_N(Z_{sN}) -\gamma(z(s)) ds \Bigr| \Bigr] 
\nonumber \\ & 
\le \sqrt{N} \Bigl( \delta_N + C_2 \Bigl\lVert \sup_{s \le 1} \left| z(s) - Z_{s}^N \right| \Bigr\rVert_2 \Bigl\lVert \frac{T^*_N}{N} - T^* \Bigr\rVert_2 \Bigr)
\nonumber \\ &
\le  C_4 \sqrt{N} \omega_N^2
\end{align}
 
After a Taylor expansion around $T^*$, it follows for the third term that
\begin{equation*}
\sqrt{N}\int_{T^*}^{T_N^*/N} \gamma(z(s)) ds 
= \sqrt{N} \gamma(1) \Bigl( \frac{T_N^*}{N} - T^* \Bigr) + \sqrt{N} R_N,
\end{equation*}
and we can subsequently bound using Proposition~\ref{prop:Bound_on_L2_norm_of_TNstar_minus_Tstar}
$$  
\sqrt{N} \E | R_N| 
\le \tfrac{1}{2} ( \sup |\gamma'| ) \E \Bigl[ \Bigl( \frac{T_N^*}{N} - T^* \Bigr)^2 \Bigr] 
\le C_5 \sqrt{N}  \omega_N^2.
$$

Using another Taylor expansion, we can expand the fourth term as
\begin{align*}
\sqrt{N}\int_0^{T^*} \gamma ( Z_s^N ) - \gamma(z(s)) ds 
&
= \sqrt{N}\int_0^{T^*} \gamma'(z(s)) ( Z_s^N - z(s) ) + \tilde R_N(s) ds 
\nonumber \\ &
= \int_0^{T^*}   \gamma'(z(s)) W_s^N ds + \sqrt{N}\int_0^{T^*} \tilde  R_N(s) ds
\end{align*}
and then bound the fourth term using Proposition~\ref{prop:LLN_comparing_Znt_and_Zt},
\begin{equation}
\sqrt{N} \E[ \tilde{R}_N ] 
\le \tfrac{1}{2} \sqrt{N} ( \sup \gamma'' ) \lVert Z_s^N - z(s) \rVert_2^2 
\le C_5 \sqrt{N} \omega_N^2.
\end{equation}

Finally, using Proposition~\ref{prop:Diffusion_approximation_with_error_bound}, we obtain for the third term:
\begin{equation}
\sqrt{N} \E[ | M^N_{T_N^*/N} - M^N_{T^*} | ]
\le \sqrt{N} \E[ < M^N_{T_N^*/N} - M^N_{T^*} > ]^{1/2} 
\le \sqrt{N} C_N^{1/2} \Bigl\lVert \frac{T_N^*}{N} - T^* \Bigr\rVert_1 
\le \sqrt{N} \omega_N^2. \nonumber
\end{equation}

Remark that $\sqrt{N} \omega_N^2 \to 0$ when $\delta_N = o(N^{-1/4})$ and $\bar{\psi}_N = o(N^{3/4})$, and then all these terms converge to $0$ as $N \to \infty$. This proves that since $\gamma(1) \neq  1$, the limit is a Gaussian random variable with variance
\begin{equation}
\sigma^2 = E\Bigl[ \Bigl( \frac{W_{T^*}}{1- \gamma(1)} \Bigr)^2 \Bigr].
\end{equation}
Defining $m_{t} = E[W_{t}^2]$, we find using It\^{o}'s formula
\begin{equation*}
\E[W_t^2]
= \E \Bigl[ \int_0^t 2 W_s d W_s + (1/2)2 \beta_{t} \Bigr]
= 2\int_0^{t} \gamma'(z(s)) E(W_s^2) ds + \beta(t)
\end{equation*}
and hence $m_t$ satisfies the differential system
\begin{gather}
\dot{m_t}=-2 \dot{ \gamma}(z_t) m_t+ \dot{\beta}(t),
\quad m_0 = 0.
\end{gather}
This finishes the proof.
\ep

\begin{example}[ER case -- Continued]
For the ER graph $\gamma(x)=\psi(x)=c(1-x)$ and $z(t) = ((1+c)/c) (1-e^{-ct}) $. Solving $z(T^*) = 1$ gives $T^* = \ln{(1+c)}/c$. The differential equation for $\beta(t)$ is given by
\begin{gather*}
\dot{\beta}(t) 
= \psi(z(s))=(1+c)e^{-ct}-1,
\end{gather*}
and the solution to \eqref{eq:ode_var_er_disc} is then
\begin{equation}
m_t 
= \exp(- 2 ct)(1-\exp(c t)  )(\exp(ct)-2c-1)\frac{1}{2c},
\end{equation}
ultimately leading to
\begin{align}
\sigma^2 
= \frac{m_{T^*}}{ ( 1-\gamma(1) )^2 } 
= \frac{c}{2(c+1)^2}.
\end{align}


\end{example}

\subsection{Continuous time version}

Our arguments that have led to our results for discrete time can be used in a similar fashion to obtain results for continuous time. We state the convergence results for fixed time intervals without proof.

\begin{prop}
Suppose that there exist a function on $\mathbb R^+$, $\gamma$ such that
\begin{gather}
\gamma_N(x)= \gamma\(\frac{x}{N}\) + \delta_N \quad \text{with} \quad \delta_N \underset{N\rightarrow \infty}{\to} 0.
\label{eq:def_gamma_cont}
\end{gather}
and suppose that the function $g(x)=(1-x)\gamma(x)$ is a ($C_L$)-Lipschitz function such that
\begin{gather}
g_N(x)= \(1-\frac{x}{N}\)\gamma_N(x)= g\(\frac{x}{N}\) + \alpha_N \quad \text{with} \quad \alpha_N \underset{N\rightarrow \infty}{\to} 0.
\label{eq:def_g}
\end{gather}
Then
\begin{gather*}
E\(\sup_{t \in [0,T]}|Z^{N}(t)- z(t)|\) \leq (\alpha'_N +2\sqrt{\lambda C_N T})\exp\{\lambda C_L T\},
\end{gather*}
where $\alpha'_N$ goes to zero with $N$, $C_N= \frac{\bar \psi_N}{N}+\frac{(1+\bar\gamma_N)^2}{N}$ and $z(t)$ is the solution of the following
(deterministic) differential equation
\begin{gather}
\dot{z}= \lambda (1-z) (1+\gamma(z)).
\label{eq:ODE}
\end{gather}
\label{prop:fluid_limit}
Suppose that there exists a function $\psi$  on $\mathbb R_+$ such that
$$\lim_{N \to \infty} \Bigl| \psi_N(x) -\psi \Bigl( \frac{x}{N} \Bigr) \Bigr| = 0.$$
The process $$W^N_t=\sqrt{N}\(Z^N_t-z(t)\)$$ converges in distribution towards $W$ which is the solution
of the following stochastic differential equation:
\begin{gather}
dW(t)=\lambda[-(1+\gamma(z(t)))+\gamma'(z(t))(1-z(t))]W(t)dt+  \sqrt{\beta'(t)} dB(t),
\end{gather}
where $B(t)$ is a standard Brownian motion, $\beta(t)=  \lambda \int_{0}^t {(1-z(s))} (\psi(z(s))+ (1+\gamma(z(s))^2 )ds$,  and $z(t)$ is the solution of \eqref{eq:ODE}.
\end{prop}

\bibliographystyle{plain}
\bibliography{probab2}

\end{document}